\documentclass[preprint,12pt]{elsarticle}

\usepackage[english]{babel}
\usepackage{hyperref}
\usepackage{vmargin}
\usepackage[ansinew]{inputenc}
\usepackage[T1]{fontenc}
\usepackage{lmodern}
\usepackage{amsfonts,amsmath,amstext,amsthm,amssymb}
\usepackage{calc}

\def\N{\mathbb{N}}
\def\Z{\mathbb{Z}}

\newtheorem{proposition}{Proposition}[section]
\newtheorem{theorem}[proposition]{Theorem}

\newtheorem{lemma}[proposition]{Lemma}
\newtheorem{corollary}[proposition]{Corollary}
\newtheorem{definition}[proposition]{Definition}

\journal{European Journal of Combinatorics}

\begin{document}

\begin{frontmatter}



\title{An Algorithmic Approach to the Extensibility of Association Schemes}


\author{M. Arora}
\address{Hausdorff Center for Mathematics, University of Bonn,\\ Endenicher Allee 62, 53115 Bonn, Germany. Email: m.arora@hcm.uni-bonn.de.}

\author{P.-H. Zieschang}
\address{Department of Mathematics, University of Texas at Brownsville,\\ Brownsville, TX 78520, USA. Email: zieschang@utb.edu.}

\address{}

\begin{abstract}
An association scheme which is associated to a height $t$ presuperscheme is said to be extensible to height $t$. Smith (1994, 2007) showed that an association scheme $\mathfrak{X}=(Q,\Gamma)$ of order $d:=|Q|$ is Schurian iff $\mathfrak{X}$ is extensible to height $(d-2)$. In this work, we formalize the maximal height $t_{\max}(\mathfrak{X})$ of an association scheme $\mathfrak{X}$ as the largest number $t\in \N$ such that $\mathfrak{X}$ is extensible to height $t$ (we also include the possibility $t_{\max}(\mathfrak{X})=\infty$, which is equivalent to $t_{\max}(\mathfrak{X})\geq (d-2)$). Intuitively, the maximal height provides a natural measure of how close an association scheme is to being Schurian.

For the purpose of computing the maximal height, we introduce the association scheme extension algorithm. On input an association scheme $\mathfrak{X}=(Q,\Gamma)$ of order $d:=|Q|$ and a number $t\in \N$ such that $1\leq t \leq (d-2)$, the association scheme extension algorithm decides in time $d^{O(t)}$ if the scheme $\mathfrak{X}$ is extensible to height $t$. In particular, if $t$ is a fixed constant, then the running time of the association scheme extension algorithm is polynomial in the order of $\mathfrak{X}$.

The association scheme extension algorithm is used to show that all non-Schurian association schemes up to order $26$ are completely inextensible, i.e. they are not extensible to any positive height $t\in \N_{>0}$. Via the tensor product of association schemes, the latter result gives rise to a multitude of examples of infinite families of completely inextensible association schemes.

\end{abstract}


\end{frontmatter}



\section{Introduction}\label{Introduction}

An association scheme $\mathfrak{X}=(Q,\Gamma)$ is composed of a finite, nonempty set $Q$ and a partition $\Gamma$ of the direct square $Q^2$, satisfying a certain set of combinatorial conditions (see \cite{BI, ZIE}). For $t\in\N$, a height $t$ presuperscheme $(Q,\Gamma^{*})$ (short: $t$-prescheme) is composed of a family of sets $\{\Gamma^{s}\}_{0\leq s \leq t}$, where each set $\Gamma^{s}$ is a partition of the direct power $Q^{s+2}$, satisfying a set of higher-dimensional variants of the combinatorial conditions which hold for association schemes (see \cite{WOJ2, WOJ, WOJ3}). Especially, if $(Q,\Gamma^{*})$ is a $t$-prescheme, then $(Q,\Gamma^0)$ constitutes an association scheme. We say that $(Q,\Gamma^0)$ is associated to $(Q,\Gamma^{*})$.

An association scheme which is associated to a $t$-prescheme (for some $t\in\N$) is said to be extensible to height $t$. Naturally, if an association scheme is extensible to height $t$, then it is also extensible to any height $t'\in \N$ with $0\leq t'\leq t$. A fundamental result by Smith \cite{SMI, SMI2} states that an association scheme $\mathfrak{X}=(Q,\Gamma)$ of order $d:=|Q|$ is Schurian iff $\mathfrak{X}$ is extensible to height $(d-2)$. A natural question arises: Given an arbitrary association scheme $\mathfrak{X}=(Q,\Gamma)$, what is the maximal height $t_{\max}(\mathfrak{X})\in \N$ which $\mathfrak{X}$ can be extended to? (Note that we include the possibility $t_{\max}(\mathfrak{X})=\infty$, which is equivalent to $t_{\max}(\mathfrak{X})\geq (d-2)$). The number $t_{\max}(\mathfrak{X})$ may provide an intuitive measure of how close the scheme $\mathfrak{X}$ is to being Schurian. 

In this paper, we introduce the association scheme extension algorithm, which on input an association scheme $\mathfrak{X}=(Q,\Gamma)$ of order $d:=|Q|$ and a number $t\in \N$ such that $1\leq t \leq (d-2)$, decides in time $d^{O(t)}$ if the scheme $\mathfrak{X}$ is extensible to height $t$. Furthermore, if $\mathfrak{X}$ is extensible to height $t$, then the algorithm outputs its unique coarsest $t$-extension $\mathfrak{X}_t$, which represents the most `basic' way in which $\mathfrak{X}$ can be extended to a $t$-prescheme. Observe that for a fixed constant $t$, the running time of the association scheme extension algorithm is polynomial in the order of $\mathfrak{X}$. 

The association scheme extension algorithm can be used to compute the maximal height $t_{\max}(\mathfrak{X})$ of a given association scheme $\mathfrak{X}$; in the worst case (i.e. if $t_{\max}(\mathfrak{X})$ is large), this may take time exponential in the order of $\mathfrak{X}$. However, heuristics suggest that association schemes of constant maximal height are far more common than association schemes of large maximal height (in a similar sense as heuristics suggest non-Schurian schemes to be far more common than Schurian schemes). Since the algorithm computes the maximal height $t_{\max}(\mathfrak{X})$ of an association scheme $\mathfrak{X}$ in polynomial time if $t_{\max}(\mathfrak{X})$ is constant, it can be quite efficient in practice.

We use the association scheme extension algorithm to prove that all non-Schurian association schemes $\mathfrak{X}=(Q,\Gamma)$ of order $|Q|\leq 26$ cannot be extended to a positive height $t\in \N_{>0}$ (schemes of the latter type are called completely inextensible). Drawing on this result, the tensor product of association schemes yields a multitude of examples of infinite families of completely inextensible association schemes.

\subsection{Related Notions}\label{Related Notions}

The focus of this work is on the notion of association schemes and their relation to $t$-preschemes. Other notions of combinatorial schemes, which include \textit{cellular algebras} (Weisfeiler \textit{et al.} \cite{WEI}), \textit{coherent configurations} (Higman \cite{HI}), \textit{Krasner algebras} (Krasner \cite{K36}), \textit{superschemes} (Smith \cite{SMI, SMI2}) and \textit{$m$-schemes} (Ivanyos, Karpinski \& Saxena \cite{IKS}) are closely related and in some cases differ only slightly, or just in notation. Moreover, the concept of $t$-preschemes and the association scheme extension algorithm seem to be closely connected to the notion of \textit{stable partitions} \cite{EP1} and the \textit{$k$-dimensional Weisfeiler-Lehman algorithm} \cite{CFI, WL68}. Note that we do not provide an exact description of the nature of these connections, or attempt to unify the various notions, as this is beyond the scope of this text.

We remark that the notion of extensibility of association schemes has gained interest in connection with recent scheme-theoretic approaches to the computational problem of factoring polynomials over finite fields \cite{AIKS, IKS}. For this line of research, it is of particular interest to gain a more thorough understanding of the combinatorial properties possessed by association schemes which are extensible to a certain height.

\subsection{Organization}\label{Organization}

\S \ref{Definitions} provides an introduction to the theory of combinatorial schemes. \S \ref{Association Schemes} defines the notion of association schemes and the concept of Schurity. \S \ref{Height $t$ Presuperschemes (short: $t$-Preschemes)} introduces the notion of $t$-preschemes and defines the concept of extensibility of association schemes. In \S \ref{Adjacency Tensors of $t$-Preschemes}, we introduce adjacency tensors of $t$-preschemes and delineate in which sense they express a central combinatorial property of $t$-preschemes (see Theorem \ref{Composition Property Theorem}). In \S \ref{The Association Scheme Extension Algorithm}, we give a description of the association scheme extension algorithm. Moreover, we list the computational results obtained through the application of the algorithm (see \S \ref{Computational Results}).

\section{Definitions}\label{Definitions}

In this section, we introduce the necessary background from the theory of combinatorial schemes to understand the context of this work. We give a survey of \textit{association schemes} (see Section \ref{Association Schemes}) and \textit{height $t$ presuperschemes} (short: \textit{$t$-preschemes}, see Section \ref{Height $t$ Presuperschemes (short: $t$-Preschemes)}). Moreover, we formalize the concept of extensibility of association schemes.

\subsection{Association Schemes}\label{Association Schemes}

\begin{definition}[Association Scheme]\label{Association Scheme Definition}
Let $Q$ be a finite nonempty set. Then an association scheme $\mathfrak{X}=(Q,\Gamma)$ on $Q$ is a partition $\Gamma=\{C_{1},...,C_{s}\}$ of the direct square $Q\times Q$, such that: \begin{enumerate}
		\item[(A1)] (Identity Relation) $C_{1}:=\{(x,x) \,|\, x\in Q\}$; 
		\item[(A2)] (Transposition) $\forall C_{i} \in \Gamma$, $C^{*}_{i}:=\{(y,x) \,|\, (x,y)\in C_{i}\} \in \Gamma$; 
	  \item[(A3)] (Intersection) $\forall C_{i}\in \Gamma$, $\forall C_{j}\in \Gamma$, $\forall C_{k}\in \Gamma$, $\exists c(i,j,k)\in \N$. $\forall (x,y)\in C_{k}$, \[
		\left|\{z\in Q \,|\, (x,z)\in C_{i},\ (z,y)\in C_{j}\}\right|=c(i,j,k).
	\] \end{enumerate} 
We refer to the numbers $c(i,j,k)$ as the intersection numbers of $\mathfrak{X}$. Moreover, we call $\left|Q\right|$ the order of $\mathfrak{X}$. \end{definition}
	
A classical example of association schemes is provided by \textit{Schurian association schemes}, which are defined below. Let $Q$ be a finite nonempty set and let $G$ be a transitive permutation group on $Q$. Let $\Gamma:=\{C_{1},...,C_{s}\}$ denote the set of orbits of $Q\times Q$ under the diagonal action of $G$, where $C_{1}:=\{(x,x) \,|\, x\in Q\}$ is the trivial orbit. Then $(Q,\Gamma)$ is an association scheme. We call schemes that arise from the action of a permutation group in the above-described manner Schurian association schemes. 

Schurian schemes provide copious examples of association schemes, but they do not cover all association schemes. A list of non-Schurian association schemes of small order can be found in Hanaki and Miyamoto's work \cite{HM}. Examples of infinite families of non-Schurian association schemes can for instance be found in \cite{EP3, FKM}.

Finding a polynomial-time algorithm which decides whether a given association scheme is Schurian or non-Schurian is a long-standing open problem. The methods introduced in \cite{BKL, BL} yield subexponential-time algorithm for testing Schurity of association schemes; this is currently the best known. Recently, Ponomarenko \cite{PON} devised an algorithm which decides the Schurity problem for \textit{antisymmetric} association schemes in polynomial time (note that an association scheme $\mathfrak{X}=(Q,\Gamma)$ is called antisymmetric if for all $C_1\neq C_{i} \in \Gamma$, $C^{*}_{i}=\{(y,x) \,|\, (x,y)\in C_{i}\}\neq C_{i}$). 

\subsection{Height $t$ Presuperschemes (short: $t$-Preschemes)}\label{Height $t$ Presuperschemes (short: $t$-Preschemes)}  

Below, we introduce the notion of \textit{height $t$ presuperschemes} (short: \textit{$t$-preschemes}), which may be regarded as a higher-dimensional analog of the notion of association schemes. In the following, let $Q$ be a finite nonempty set. For each $n\in \N_{>1}$, define a projection \begin{align*}
	\text{pr}_n: Q^n &\longrightarrow Q^{n-1}\\
	(x_1,...,x_{n-1},x_n)&\longrightarrow (x_1,...,x_{n-1})
\end{align*} (the projection $\text{pr}_n$ eliminates the last coordinate from tuples in $Q^n$). The inverse image of a set $C\subseteq Q^{n-1}$ under $\text{pr}_n$ is denoted by $\text{pr}^{-1}_n (C)$. Throughout this work, we omit the index $n$ (we assume it is clear from context) and just write $\text{pr}$ instead of $\text{pr}_n$. For each $n\in \N$, observe that the symmetric group on $n$ elements $\text{Symm}_n$ acts on the set of tuples $Q^n$ by permuting the coordinates. For all $\bar{u}:=(u_1,...,u_n)\in Q^n$ and $\tau\in \text{Symm}_n$, define \[
	\bar{u}^\tau:=(u_{\tau(1)},...,u_{\tau(n)}).	
\] Furthermore, we fix the following convention: \[
  \N_{t}:=\{n\in \N \,|\, n\leq t\}, \ \ \ \N^{2}_{t}:=\{(m,n)\in \N^{2} \,|\, m+n\leq t\}. 
 \] Note that the definition of height $t$ presuperschemes given below is equivalent to the definition given by Wojdy\l\text{o} \cite{WOJ2, WOJ, WOJ3}.
 
 \begin{definition}[Height $t$ Presuperscheme]\label{Presuperscheme Definition}
Let $Q$ be a finite nonempty set and let $t\in \N$. A height $t$ presuperscheme $(Q,\Gamma^{*})$ on $Q$ is a family of sets $\{\Gamma^{n}\}_{n\in \N_t}$, where each set $\Gamma^{n}=\{C^{n}_{1},...,C^{n}_{s_{n}}\}$ is a partition of the direct power $Q^{n+2}$ (note that all $C^{n}_{i}$ are assumed to be nonempty), such that: \begin{enumerate}

   \item[(P1)] (Identity Relation) $C^{0}_{1}:=\{(x,x) \,|\, x\in Q\}$; 

   \item[(P2)] (Projection) $\forall n\in \N_{t}-\{0\}$, $\forall C^{n}_{j} \in \Gamma^{n}$, \[
  \mathrm{pr} (C^{n}_{j}):=\{\mathrm{pr}(\bar{u}) \,|\, \bar{u}\in C^{n}_{j}\}\in \Gamma^{n-1};
 \] 

   \item[(P3)] (Invariance) $\forall n\in \N_{t}$, $\forall C^{n}_{j} \in \Gamma^{n}$, $\forall \tau\in \mathrm{Symm}_{n+2}$, \[
  (C^{n}_{j})^\tau:=\{\bar{u}^\tau \,|\, \bar{u}\in C^{n}_{j}\}\in \Gamma^n;
 \]

   \item[(P4)] (Intersection) $\forall (m,n)\in \N^{2}_{t}$, $\forall C^{m}_{i}\in \Gamma^{m}$, $\forall C^{n}_{j}\in \Gamma^{n}$, $\forall C^{m+n}_{k}\in \Gamma^{m+n}$, \linebreak $\exists c(i,j,k;m,n)\in \N$. $\forall (x_{0},...,x_{m},y_{0},...,y_{n})\in C^{m+n}_{k}$, \[
  \left|\{z\in Q \,|\, (x_{0},...,x_{m},z)\in C^{m}_{i},\ (z,y_{0},...,y_{m})\in C^{n}_{j}\}\right|=c(i,j,k;m,n).
\] 

\end{enumerate} For brevity, we refer to height $t$ presuperschemes simply as \textit{$t$-preschemes}. We call the elements of $\Gamma^{n}$ ($0\leq n \leq t$) the \textit{relations at height $n$}. We refer to the numbers $c(i,j,k;m,n)$ as the \textit{intersection numbers} of $(Q,\Gamma^{*})$. \end{definition} 

Property (P2) interrelates the different layers $\{\Gamma^{n}\}_{n\in \N_t}$ of a $t$-prescheme, while \linebreak Properties (P3), (P4) may be regarded as higher-dimensional analogs of Properties (A2), (A3) of association schemes, respectively. From Definition \ref{Presuperscheme Definition} it is clear that a $0$-prescheme and an association scheme constitute the exact same notion. 

If $(Q,\Gamma^{*})$ is a $t$-prescheme, then $(Q,\Gamma^{0})$ is an association scheme. We say that the association scheme $(Q,\Gamma^{0})$ is \textit{associated} to the $t$-prescheme $(Q,\Gamma^{*})$. If an association scheme $\mathfrak{X}$ is associated to a $t$-prescheme $(Q,\Gamma^{*})$, we call $\mathfrak{X}$ \textit{extensible to height $t$}. In this case, we refer to the $t$-prescheme partitions $\{\Gamma^{n}\}_{1\leq n \leq t}$ as a \textit{$t$-extension} of $\mathfrak{X}$. Note that by definition, every association scheme is extensible to height $0$. 

We define the \textit{maximal height} $t_{\max}(\mathfrak{X})$ of an association scheme $\mathfrak{X}$ as the largest number $t\in \N$ such that $\mathfrak{X}$ is extensible to height $t$. If $\mathfrak{X}$ is extensible to arbitrary heights (meaning that for all $t\in \N$, $\mathfrak{X}$ is extensible to height $t$), we say that $\mathfrak{X}$ has maximal height $\infty$. In case $t_{\max}(\mathfrak{X})=0$, we say that $\mathfrak{X}$ is \textit{completely inextensible}.

For an association scheme $\mathfrak{X}=(Q,\Gamma)$ of order $d:=|Q|$, it is easily proven that $t_{\max}(\mathfrak{X})=\infty$ iff $\mathfrak{X}$ is extensible to height $(d-2)$. A fundamental result by Smith connects the concept of extensibility to the notion of Schurity of association schemes.

\begin{theorem}[Smith \cite{SMI, SMI2}]\label{Smith Result}
An association scheme $\mathfrak{X}=(Q,\Gamma)$ of order $d:=|Q|$ is Schurian iff $\mathfrak{X}$ is extensible to height $(d-2)$.
\end{theorem} 

Note that Theorem \ref{Smith Result} may also be phrased as follows: An association scheme $\mathfrak {X}$ is Schurian iff $t_{\max}(\mathfrak{X})=\infty$. Moreover, observe that if an association scheme $\mathfrak {X}=(Q,\Gamma)$ of order $d:=|Q|$ is non-Schurian, then $0 \leq t_{\max}(\mathfrak{X}) < (d-2)$.

\section{Adjacency Tensors}\label{Adjacency Tensors of $t$-Preschemes}
  
In this section, we introduce the notion of \textit{adjacency tensors}. The concept of adjacency tensors of $t$-preschemes naturally generalizes the notion of \textit{adjacency matrices} of association schemes (see \cite{BI, ZIE}). Analogously, adjacency tensors describe the intersection property of $t$-preschemes in simple algebraic terms (see Theorem \ref{Composition Property Theorem}). \linebreak We apply the notion of adjacency tensors in Section \ref{The Association Scheme Extension Algorithm}, when we describe the association scheme extension algorithm.

\subsection{$k$-Tensors}\label{$k$-Tensors and Tensor Composition}

In the following, we introduce \textit{tensors of order $k$} (short: \textit{$k$-tensors}) and define certain natural operations associated with this notion. Note that $k$-tensors constitute a natural generalization of the concept of square matrices.

\begin{definition}[$k$-Tensor]\label{Rank $k$ Tensor}
For $k\geq 2$, a \textit{$k$-tensor} with entries in $\Z$ is a function \[
   T:\{1,...,d\}^{k}\longrightarrow \Z. 
	\] We refer to the number $k$ as the order of the tensor $T$. We denote by $T_{i_{1}\cdots i_{k}}$ the image of $(i_{1},...,i_{k})$ under $T$. We call $T_{i_{1}\cdots i_{k}}$ the \textit{$(i_{1},...,i_{k})$-entry} of $T$.
\end{definition}

In this work, tensors are regarded simply as multidimensional arrays. For $k=2$, the notion of $k$-tensors with entries in $\Z$ coincides with the notion of $d\times d$ matrices with entries in $\Z$. For a more general (algebraic) treatment of tensors, the reader is referred to \cite{CL, DIM}.

In the following, we define some basic operations for $k$-tensors. These operations naturally generalize the standard matrix operations from linear algebra. First, for two $k$-tensors $S,T:\{1,...,d\}^{k}\longrightarrow \Z$, we define their \textit{sum} $U=S+T$ as the $k$-tensor $U:\{1,...,d\}^{k}\longrightarrow \Z$ with entries \[
   U_{i_{1}\cdots i_{k}}= S_{i_{1}\cdots i_{k}} + T_{i_{1}\cdots i_{k}}.
   \] Next, for an element $c\in \Z$ and a $k$-tensor $S:\{1,...,d\}^{k}\longrightarrow \Z$, we define their \textit{scalar product} $V=c\cdot S$ as the $k$-tensor $V:\{1,...,d\}^{k}\longrightarrow \Z$ with entries \[
   V_{i_{1}\cdots i_{k}}= c \cdot S_{i_{1}\cdots i_{k}}.
   \] Moreover, for a $m$-tensor $E:\{1,...,d\}^{m}\longrightarrow \Z$ and a $n$-tensor $F:\{1,...,d\}^{n}\longrightarrow \Z$, \linebreak we define their \textit{inner product} $W=EF$ as the rank $(m+n-2)$ tensor \linebreak $W:\{1,...,d\}^{(m+n-2)}\longrightarrow \Z$ with entries \[
   W_{i_{1}\cdots i_{m+n-2}}=\sum^{d}_{j=1} E_{i_{1}\cdots i_{m-1}j} \cdot F_{ji_{m}\cdots i_{m+n-2}}.
   \] The above operations generalize the standard addition, scalar multiplication and inner multiplication of matrices. It is easily verified that addition and inner multiplication of tensors are associative, distributive and compatible with scalar multiplication. \linebreak

\vspace{-0.5cm}

\subsection{Adjacency Tensors of $t$-Preschemes}\label{Adjacency Tensors}

In the following, we define the notion of adjacency tensors, boolean tensors which indicate membership to subsets of direct powers of $Q:=\{1,...,d\}$.

\begin{definition}[Adjancency Tensor]\label{Adjacency Tensor}
Let $Q:=\{1,...,d\}$ and let $C\subseteq Q^n$, where $n\geq 2$. We define the adjacency tensor corresponding to the subset $C$ as the $n$-tensor $A(C):\{1,...,d\}^{n}\longrightarrow \Z$  such that the component $[A(C)]_{x_{1}\cdots x_{n}}$ is $1$ if $(x_{1},...,x_{n})\in C$ and $0$ otherwise.  
\end{definition}

Let $(Q,\Gamma^{*})$ be a $t$-prescheme on $Q:=\{1,...,d\}$. We denote the \textit{adjacency tensor of a relation} $C^{m}_{i}\in \Gamma^{m}$ ($m\in \N_t$) as the $(m+2)$-tensor $A^{m}_{i}:\{1,...,d\}^{m+2}\longrightarrow \Z$, where $(A^{m}_{i})_{x_{1}\cdots x_{m+2}}$ is $1$ if $(x_{1},...,x_{m+2})\in C^{m}_{i}$ and $0$ otherwise. Adjacency tensors can be used to express the intersection property of $t$-preschemes in algebraic terms (analogously to \textit{adjacency matrices} in the case of association schemes, see \cite{BI, ZIE}).

\begin{theorem}\label{Composition Property Theorem}
Let $(Q,\Gamma^{*})$ be a $t$-prescheme on the set $Q:=\{1,...,d\}$. Then for all $(m,n)\in \N^{2}_{t}$, $C^{m}_{i}\in \Gamma^{m}$ and $C^{n}_{j}\in \Gamma^{n}$, it holds that \[
   A^{m}_{i}A^{n}_{j}=\sum^{s_{m+n}}_{k=1} c(i,j,k;m,n)A^{m+n}_{k},
   \] where $A^{m}_{i},A^{n}_{j}$ and $A^{m+n}_{k}$ denote the adjacency tensors of $C^{m}_{i},C^{n}_{j}$ and $C^{m+n}_{k}\in \Gamma^{m+n}$,\linebreak respectively, and $c(i,j,k;m,n)\in \N$ denote the intersection numbers. Moreover, \linebreak the above statement is equivalent to the intersection property of $t$-preschemes \linebreak (see Definition \ref{Presuperscheme Definition} (P4)).
\proof
The intersection property of $t$-preschemes states that for all $(m,n)\in \N^{2}_{t}$, \linebreak $C^{m}_{i}\in \Gamma^{m}$, $C^{n}_{j}\in \Gamma^{n}$, $C^{m+n}_{k}\in \Gamma^{m+n}$ and $(x_{0},...,x_{m},y_{0},...,y_{m})\in C^{m+n}_{k}$, it holds that \[
	 c(i,j,k;m,n)=\left|\{z\in Q \,|\, (x_{0},...,x_{m},z)\in C^{m}_{i},\ (z,y_{0},...,y_{m})\in C^{n}_{j}\}\right|.
	\] Note that the above equation can also be written as \[
   c(i,j,k;m,n)=\sum^{d}_{z=1}\left(A^{m}_{i}\right)_{x_{0} \cdots x_{m}z}\left(A^{n}_{j}\right)_{zy_{0} \cdots y_{m}}
   \] where the right-hand side is $\left(A^{m}_{i}A^{n}_{j}\right)_{x_{0} \cdots x_{m}y_{0} \cdots y_{m}}$ by the definition of the inner product of tensors. From this the assertion follows immediately. \qed 
\end{theorem}

\section{The Association Scheme Extension Algorithm}\label{The Association Scheme Extension Algorithm}

In this section, we describe the association scheme extension algorithm. On input an association scheme $\mathfrak{X}=(Q,\Gamma)$ of order $d:=|Q|$ and a number $t\in \N$ such that \linebreak $1\leq t \leq (d-2)$, the association scheme extension algorithm decides in time $d^{O(t)}$ if $\mathfrak{X}$ is extensible to height $t$. Furthermore, if $\mathfrak{X}$ is extensible to height $t$, then the algorithm outputs its unique coarsest $t$-extension $\mathfrak{X}_t$, which represents the most `basic' way in which $\mathfrak{X}$ can be extended to a $t$-prescheme. We apply the association scheme extension algorithm to determine that all non-Schurian association schemes up to order $26$ are completely inextensible (see Theorem \ref{Comp In}). Via the tensor product of association schemes, the latter result gives rise to a multitude of examples of infinite families of completely inextensible association schemes (see Section \ref{Computational Results}).

\subsection{Description of the Algorithm}\label{Description of the Algorithm}

We now describe the association scheme extension algorithm. On input an association scheme $\mathfrak{X}=(Q,\Gamma)$ on $Q:=\{1,...,d\}$ and a number $t\in \N$ such that $1\leq t \leq (d-2)$, the algorithm begins with trivial partitions $\Gamma^s:=\{Q^{s+2}\}$ ($1\leq s\leq t$) and then gradually refines these partitions according to a set of rules derived from the properties of $t$-extensions (see Definition \ref{Presuperscheme Definition}). Via this refinement process, the partitions $\Gamma^s$ ($1\leq s\leq t$) either turn into a $t$-extension of $\mathfrak{X}$, or they provide combinatorial justification for the conclusion that $\mathfrak{X}$ cannot be extended to height $t$.

\vspace{0.25cm}

\noindent \textbf{Input:} An association scheme $\mathfrak{X}=(Q,\Gamma)$ on $Q:=\{1,...,d\}$, and a number $t\in \N$ such that $1\leq t \leq (d-2)$.

\vspace{0.25cm}

\noindent \textbf{Output:} A $t$-extension $\{\Gamma^{s}\}_{1\leq s \leq t}$ of $\mathfrak{X}$, or the decision that $\mathfrak{X}$ is not extensible to height $t$.

\vspace{0.25cm}

\noindent \textbf{Initialization.} For each $1\leq s\leq t$, let $\Gamma^s:=\{Q^{s+2}\}$ be the trivial partition of $Q^{s+2}$. \linebreak

\vspace{-0.25cm}

\noindent \textbf{Step 1.} For each $1\leq s\leq t$, refine the partition $\Gamma^s$ of $Q^{s+2}$ according to the pro-jection property of $t$-preschemes (see Definition \ref{Presuperscheme Definition} \textit{(P2)}). That is, for each $C\in \Gamma^s$, determine if the set $\text{pr}(C)$ can be written as a union of relations in $\Gamma^{s-1}$, i.e. if \[
  \text{pr}(C)=C^{s-1}_{i_{1}}\cup \cdots \cup C^{s-1}_{i_{k}} 
 \] for some $C^{s-1}_{i_{1}},...,C^{s-1}_{i_{k}}\in \Gamma^{s-1}$.

\begin{enumerate}

\itemindent30pt 

   \item[If YES.] Replace in $\Gamma^s$ the set $C\in \Gamma^s$ with the pairwise disjoint sets \[
  C\cap \text{pr}^{-1}(C^{s-1}_{i_{1}}),...,C\cap \text{pr}^{-1}(C^{s-1}_{i_{k}}).
 \] 

   \item[ELSE.] Distinguish between the following two cases:

\vspace{0.2cm}

   \begin{enumerate}

\itemindent48pt 

     \item[(a) If $s>1$.] Replace in $\Gamma^{s-1}$ each set $C'\in \Gamma^{s-1}$ such that $C'\cap \text{pr}(C) \neq \emptyset$ with the two disjoint sets $C'\cap \text{pr}(C)$ and $C'\setminus \text{pr}(C)$.

\vspace{0.2cm}

     \item[(b) If $s=1$.] Terminate the algorithm and output: \textit{$\mathfrak{X}$ is not extensible to height $t$}.
   \end{enumerate}

\end{enumerate}

\vspace{0.05cm}

\noindent \textbf{Step 2.} For each $1\leq s\leq t$, refine the partition $\Gamma^s$ of $Q^{s+2}$ according to the invariance property of $t$-preschemes (see Definition \ref{Presuperscheme Definition} \textit{(P3)}). That is, for each $C\in \Gamma^s$ and each $\tau \in Symm_{s+2}$, replace in $\Gamma^s$ each set $C'\in \Gamma^s$ such that $C'\cap C^\tau \neq \emptyset$ with the two disjoint sets $C'\cap C^\tau$ and $C'\setminus C^\tau$. 

\vspace{0.25cm}

\noindent \textbf{Step 3.} For each $1\leq s\leq t$, refine the partition $\Gamma^s$ of $Q^{s+2}$ according to \linebreak the intersection property of $t$-preschemes (see Theorem \ref{Composition Property Theorem}). That is, for each $m,n\in\N$ such that $s=(m+n)$, and each pair of sets $C^{m}_{i}\in \Gamma^{m}$ and $C^{n}_{j}\in \Gamma^{n}$, compute the inner product \[
  P:=A^{m}_{i}A^{n}_{j},
\] where $A^{m}_{i},A^{n}_{j}$ denote the adjacency tensors of $C^{m}_{i},C^{n}_{j}$, respectively (see Section \ref{Adjacency Tensors of $t$-Preschemes}). The entries of $P$ are integers in the range from $0$ to $d$. For each $r=0,...,d$ define \[
  P^{-1}(r):=\{(i_1,...,i_{s+2}) \in Q^{s+2} \,|\, P_{i_1 \cdots i_{s+2}}=r\}
 \] and replace in $\Gamma^s$ each set $C\in \Gamma^s$  such that $C\cap (P^{-1}(r)) \neq \emptyset$ with the two disjoint sets $C\cap (P^{-1}(r))$ and $C\setminus (P^{-1}(r))$.

\vspace{0.25cm}

\noindent \textit{Repeat Steps 1-3. If none of them yields any further refinement of the partitions \linebreak $\Gamma^s$ ($1\leq s\leq t$), then terminate the algorithm and output $\{\Gamma^{s}\}_{1\leq s \leq t}$.} \qed

\subsection{Correctness of the Algorithm}\label{Correctness of the Algorithm}
 
We now prove the correctness of the association scheme extension algorithm. We start with a preliminary lemma.

\begin{lemma}\label{Correctness Lemma}
Let $\mathfrak{X}=(Q,\Gamma)$ be an association scheme on $Q:=\{1,...,d\}$ and let $t\in \N$ be such that $1\leq t \leq (d-2)$. The following holds: \begin{enumerate}
   \item[(1)] On input $\mathfrak{X}$ and $t$, the association scheme extension algorithm terminates after at most $d^{O(t)}$ steps.

   \item[(2)] On input $\mathfrak{X}$ and $t$, if the association scheme extension algorithm outputs a set of partitions $\{\Gamma^{s}\}_{1\leq s \leq t}$, then these partitions constitute a $t$-extension of $\mathfrak{X}$.
\end{enumerate} 

\proof
(1) Note that the algorithm can make at most $(d^3+...+d^{t+2})$ refinements to the partitions $\{\Gamma^s\}_{1\leq s\leq t}$ before it must terminate. Moreover, observe that the algorithm goes through at most $d^{O(t)}$ elementary operations in between two refinements. From this the assertion follows immediately.

(2) Note that the algorithm outputs a set of partitions $\{\Gamma^s\}_{1\leq s\leq t}$ only if \linebreak Steps 1-3 of the algorithm do not yield any further refinement of $\{\Gamma^s\}_{1\leq s\leq t}$. The latter condition implies that Definition \ref{Presuperscheme Definition} \textit{(P2)}-\textit{(P4)} hold for $\mathfrak{X}$ and $\{\Gamma^s\}_{1\leq s\leq t}$ \linebreak (see Theorem \ref{Composition Property Theorem}). This in turn implies that the partitions $\{\Gamma^s\}_{1\leq s\leq t}$ constitute a $t$-extension of $\mathfrak{X}$.
\qed 
\end{lemma}

Let us fix some terminology. Let $X$ be a finite, nonempty set and let $\mathcal{P}, \mathcal{R}$ be partitions of $X$. If for each $P\in \mathcal{P}$ there exist sets $R_1,...,R_n\in \mathcal{R}$ such that $P=\cup^{n}_{i=1} R_i$, then we call $\mathcal{P}$ a \textit{fusion} of $\mathcal{R}$. We use this convention in the proof of correctness of the association scheme extension algorithm given below.

\begin{theorem}\label{Correctness of ASE}
The association scheme extension algorithm works correctly. The running time of the association scheme extension algorithm is $d^{O(t)}$. 
\proof
Let $\mathfrak{X}=(Q,\Gamma)$ be an association scheme on $Q:=\{1,...,d\}$ and let $t\in \N$ be such that $1\leq t \leq (d-2)$. First, assume $\mathfrak{X}$ is not extensible to height $t$. Then by Lemma \ref{Correctness Lemma} (1), (2) it follows that on input $\mathfrak{X}$ and $t$, the algorithm correctly outputs the decision that $\mathfrak{X}$ is not extensible to height $t$, in time $d^{O(t)}$.

Now consider the converse: Assume we are given as input an association scheme $\mathfrak{X}=(Q,\Gamma)$ on $Q:=\{1,...,d\}$ and a number $t\in \N$ with $1\leq t \leq (d-2)$ such that $\mathfrak{X}$ is extensible to height $t$. Choose an arbitrary $t$-extension $\{\tilde{\Gamma}^s\}_{1\leq s\leq t}$ of $\mathfrak{X}$. Observe the following facts about the partitions $\{\Gamma^s\}_{1\leq s\leq t}$ which appear in the algorithm:

\begin{enumerate}
   \item[(i)] For each $1\leq s\leq t$, the partition $\Gamma^s$ is trivially a fusion of $\tilde{\Gamma}^s$ at the initialization step.

\vspace{0.1cm}

   \item[(ii)] For each $1\leq s\leq t$, the partition $\Gamma^s$ remains a fusion of $\tilde{\Gamma}^s$ over the whole course of the algorithm (this follows from Properties \textit{(P2), (P3), (P4)} of \linebreak Definition \ref{Presuperscheme Definition} applied on $\mathfrak{X}$ and $\{\tilde{\Gamma}^s\}_{1\leq s\leq t}$). Especially, the algorithm never terminates during the execution of Step 1.

\end{enumerate} 

By statement (ii) and Lemma \ref{Correctness Lemma} (1), we conclude that on input $\mathfrak{X}$ and $t$, the algorithm outputs a set of partitions $\{\Gamma^s\}_{1\leq s\leq t}$. By Lemma \ref{Correctness Lemma} (2), the output $\{\Gamma^s\}_{1\leq s\leq t}$ constitutes a $t$-extension of $\mathfrak{X}$. \qed 
\end{theorem}

Recall that in the proof of Theorem \ref{Correctness of ASE}, the $t$-extension $\{\tilde{\Gamma}^s\}_{1\leq s\leq t}$ of $\mathfrak{X}$ was chosen arbitrarily. Hence we obtain the following corollary:

\begin{corollary}\label{Unique Coarsest Extension}
On input an association scheme $\mathfrak{X}=(Q,\Gamma)$ and a number $t\in \N$ with $1\leq t \leq (d-2)$ such that $\mathfrak{X}$ is extensible to height $t$, the association scheme extension algorithm outputs the \textit {unique coarsest $t$-extension} $\mathfrak{X}_t:=\{\Gamma^s\}_{1\leq s\leq t}$ of $\mathfrak{X}$. That is, for any $t$-extension $\{\tilde{\Gamma}^s\}_{1\leq s\leq t}$ of $\mathfrak{X}$, for each $1\leq s\leq t$, the partition $\Gamma^s$ is a fusion of $\tilde{\Gamma}^s$. 
\end{corollary}

\subsection{Computational Results}\label{Computational Results}

We employ the association scheme extension algorithm to determine the extensibility properties of all non-Schurian association schemes up to order $26$. Note that there are exactly 142 non-Schurian schemes of order less or equal to $26$ (see \cite{HM2}).

\begin{theorem}\label{Comp In}
All non-Schurian association schemes $\mathfrak{X}=(Q,\Gamma)$ of order $|Q|\leq 26$ are completely inextensible. 
\proof
We created a program of the association scheme extension algorithm with fixed parameter $t=1$ in the input, written in ``C''. We applied our program to all non-Schurian association schemes of order less or equal to $26$; for this we relied on the classification of non-Schurian association schemes of small order by Hanaki and Miyamoto \cite{HM2, HM3, HM, HM4}. The reader can download an organized version of the C-programs and their output online \cite{ARO}. \qed
\end{theorem}

Let us fix some convention.  For an association scheme $\mathfrak{X}=(Q,\Gamma)$, we denote the equivalence relation on $Q\times Q$ corresponding to the partition $\Gamma$ by $\equiv _{\Gamma}$. Recall the definition of the \textit{tensor product} of association schemes. For two association schemes $\mathfrak{X}_1=(Q_1,\Gamma_1)$ and $\mathfrak{X}_2=(Q_2,\Gamma_2)$, the tensor product $\mathfrak{X}_1 \otimes \mathfrak{X}_2$ is defined as the association scheme $(Q_1 \times Q_2, \Gamma_1\otimes \Gamma_2)$ such that for all $x_1,x'_1,y_1,y'_1\in Q_1$ and $x_2,x'_2,y_2,y'_2\in Q_2$, \begin{align*} 
	&((x_1,x_2), (x'_1, x'_2))\equiv _{\Gamma_1\otimes \Gamma_2} ((y_1,y_2), (y'_1, y'_2))\\
\Longleftrightarrow &(x_1,x'_1)\equiv _{\Gamma_1} (y_1,y'_1) \ \text{and} \ (x_2,x'_2)\equiv _{\Gamma_2} (y_2,y'_2). 
\end{align*} Given a number $t\in \N$, it is easily seen that the tensor product $\mathfrak{X}_1 \otimes \mathfrak{X}_2$ is extensible to height $t$ iff both $\mathfrak{X}_1$ and $\mathfrak{X}_2$ are extensible to height $t$. Via the above construction, Theorem \ref{Comp In} gives rise to a multitude of examples of infinite families of completely inextensible association schemes. Especially, we have the following corollary.

\begin{corollary}\label{Infinitely Comp In}
There exist infinitely many completely inextensible association \linebreak schemes.
\end{corollary}

\section{Conclusion and Open Questions}\label{Conclusion and Open Questions}

For an association scheme $\mathfrak{X}=(Q,\Gamma)$, we defined the notion of extensibility to height $t$ (for a given $t\in \N$), the notion of the maximal height $t_{\max}(\mathfrak{X})$ and - assuming that $\mathfrak{X}$ is extensible to height $t$ - the concept of the unique coarsest $t$-extension $\mathfrak{X}_t$. We delineated in which sense the maximal height may be regarded as an intuitive measure of how close an association scheme is to being Schurian. Furthermore, we described the association scheme extension algorithm, which on input an association scheme $\mathfrak{X}=(Q,\Gamma)$ of order $d:=|Q|$ and a number $t\in \N$ such that $1\leq t \leq (d-2)$, decides in time $d^{O(t)}$ if $\mathfrak{X}$ is extensible to height $t$. We used the association scheme extension algorithm to determine that all non-Schurian association schemes up to order $26$ are completely inextensible, i.e. they have maximal height $0$. 

Computing the maximal height of an association scheme $\mathfrak{X}=(Q,\Gamma)$ with the association scheme extension algorithm may require time exponential in $|Q|$ in the worst case. A central open question is whether there exists an algorithm for computing the maximal height which achieves a better worst-case running time (for instance, \linebreak a running time in the subexponential range). A relaxation of this question would be to ask whether there exist `thresholds' $t(d)\in \N$ such that for all association schemes $\mathfrak{X}=(Q,\Gamma)$ of order $d:=|Q|$, deciding if $\mathfrak{X}$ is extensible to height $t(d)$ can be done more efficiently than using the association scheme extension algorithm. Apart from this, we remark that it is currently an open problem to identify the smallest order $d\in \N$ for which there exists a non-Schurian association scheme of positive maximal height. We leave the above questions to future research. 

\section*{Acknowledgments}

The authors would like to thank the Hausdorff Center for Mathematics, University of Bonn and the University of Texas at Brownsville, for their kind support and for research funding. Moreover, the authors would like to thank Marek Karpinski, Ilya Ponomarenko and Nitin Saxena for the many fruitful discussions and helpful remarks. 





\bibliographystyle{model1b-num-names}
\bibliography{refs}

\end{document}